# To substantiate the asymptotics of solving the Cauchy problem for a singularly perturbed weakly nonlinear transport equation[*]

**Nesterov A.V.** [1][0000-0002-4702-4777] and **Zaborsciy A.V.**[2][0000-0003-2480-7215]

[1] PLEKHANOV Russian University of Economics, Stremyanny lane 36, Moscow, 117997, Russia
andrenesterov@yandex.ru
[2] RPE «RADICO» Ltd , Marx Ave., 14a, Obninsk, Kaluga region, 249035, Russia
alexander.zaborskiy@mail.ru

**Abstract.** A theorem is proved on the uniform estimation of the residual term of the asymptotic expansion with respect to a small parameter of the solution of the initial problem for a singularly perturbed differential operator weakly nonlinear transport equation in the critical case.

## 1 Introduction

In the work of the authors [1] an asymptotic expansion (AE) in a small parameter was constructed for the solution of an initial problem for a singularly perturbed differential-operator equation with variable coefficients and small nonlinearity

$$\varepsilon^2(U(x,t,p)_t + D(p)U(x,t,p)_x) = L_p U(x,t,p) + \varepsilon^2 F(x,p,U) \quad (1)$$

$$U(x,0,p) = w(x/\varepsilon, p), \quad (2)$$

where:
$U(x,t,p)$ - solution, $\{x,t,p\} \in H = \{|x| < \infty; 0 \le t \le T, T > 0; p_1 \le p \le p_2\}$; $0 < \varepsilon \ll 1$ is a small parameter, $D(p)$ - a continuous of $p$ in $p \in P$ (for example $P = [p_1, p_2]$) function, satisfying the condition $|D(p)| \ge D_0 > 0$, $L_p : A_p \to A_p$ linear operator acting in the space $A_p$ of continuous on $p$ functions $U(x,t,p)$ with inner product $(h_1, h_2)$, $F(x,p,U)$ - a fairly smooth function.

[*] This research was performed in the framework of the state task in the field of scientific activity of the Ministry of Science and Higher Education of the Russian Federation, project "Development of the methodology and a software platform for the construction of digital twins, intellectual analysis and forecast of complex economic systems", grant no. FSSW-2020-0008.

$w(z,p)$ -initial condition - together with its derivatives with respect to $z$ up to the order $N+3$ satisfies the inequality $w(z,p) \leq Ce^{-\beta\|z\|^2}$, where $N$ is a natural number, $C \geq 0, \beta > 0$.

I. The operator $L_p$ has a single zero eigenvalue $\lambda_0 = 0$, $h_0(p)$ eigenfunction of the operator $L_p$ corresponding zero eigenvalue, $h_0^*(p)$ - eigenfunction of the adjoint operator $L_p^*$, corresponding zero eigenvalue $\lambda_0^* = 0$.

In [1], if a number of conditions are met (see below), the solution AE of problem (1)-(2) is constructed in the form

$$U(x,t,p) = S(\zeta,t,p) + P(\xi,\tau,p) + R =$$
$$= \sum_{i=0}^{N} \varepsilon^i (s_i \zeta, t, p) + p_i(\xi, \tau, p)) + R_N = U_N + R_N \quad (3)$$

где $S$ - surge function, $P$ - boundary function, $R$ – residual term;

$$\zeta = \frac{t - Bx}{\varepsilon}, \quad \xi = \frac{x}{\varepsilon}, \quad \tau = \frac{t}{\varepsilon^2}, \quad B = \frac{(h_0, h_0^*)}{(D(s,p)h_0, h_0^*)}.$$

In the same place, problems were obtained by which all the members of the AE (3) are determined, estimates of the members of the AE are obtained, and also, under a number of conditions, a theorem on the estimate of the remainder by the residual is proved. However, in this work, there was no estimate of the remainder in a uniform (or any other) norm. The present work makes up for this deficiency.

## 2 Estimation of the residual term $R$

We present a theorem on estimating the residual term in AE (3) by the residual [1].
Let the conditions [1]

I. The eigenvalues of the operator the operator $L_p$ have, except $\lambda_0$, negative real parts: $\text{Re}\,\lambda \leq -k$, $k > 0$.

II. $(h_0, h_0^*) \neq 0$. It is assumed below that $(h_0, h_0^*) = 1$

III. The function $g(t) = (\Psi_0 G \Psi_0 h_{00}, h_{00}^*) < 0$. Here $\Psi_0(p) = D_0(p) - (D_0(p)h_0(p), h_0^*(p))$, $G$ the pseudo-inverse to $L_p$ operator.

*Remark.* We call an operator $G$ pseudo-inverse to an operator $L_p$ if the solution of the equation $L_p U = F$ under the condition $(F, h_0^*) = 0$ can be written in the form $U = GF + Ch_0$, where $C$ does not depend on $p$.

IV. The set of eigenvalues $L_p: \lambda_0, \lambda_1, \lambda_2, \ldots$ is countable.

V. The eigenvalues of the operator $L_p : \lambda_0, \lambda_1, \lambda_2, \ldots$ - are one-time, and the corresponding eigenfunctions: $h_0, h_1, h_2, \ldots$ - are orthogonal, normalized, form a complete system of functions and satisfy the conditions $\forall i : |h_i(p)| < const$.

VI. The function $w(\xi, p)$ decomposes into an absolutely and uniformly convergent series $w(\xi, p) = \sum_{i=0}^{\infty} w_i(\xi) h_i(p)$, which can be purposefully differentiated a sufficient number of times.

Having chosen a certain natural number $N$, we define the term $N$ of AE (3) and represent the solution of problem (1) - (2) in the form (3). In [1], the theorem on the estimation of the residual term by the residual was proved

*Theorem.* Let $U(x, t, p)$ be the solution of the initial problem

$$\varepsilon^2 (U_t + D(p) U_x) = L_p U + \varepsilon^2 F(x, p, U),$$
$$U(x, 0, p) = w(x/\varepsilon, p),$$

Let conditions I-VI be met. Then the solution of problem (1)-(2) can be represented in the form $U = U_N + R_N$, where $U_N$ is the above-constructed AR of the solution, and the residual term $R$ satisfies the asymptotic estimate by the residual:

$$\varepsilon^2 (R_{Nt} + D R_{Nx}) = L_p R_N + \varepsilon^2 RF + O(\varepsilon^{N+2}),$$
$$R_N |_{t=0} = 0. \tag{4}$$

This theorem, however, does not provide a rigorous justification for the constructed method for solving problem (1)-(2), since it is not clear from problem (4) whether there is an estimate of the value itself $R$ in the uniform (or any other) norm.

This paper fills this gap.

We estimate the norm of the residual term $R_{N+3}$ under additional conditions VII, VIII:

*VII.* There is a constant $K$ that for any positive function $u(p) > 0$, the value $L_p u - Ku$ is also positive $L_p u - Ku > 0$.

*VIII.* The eigenfunction $h$ corresponding to the eigenvalue $\lambda = 0$ is continuous with respect to $p, p \subset [p_1, p_2]$ and sign-constant (and, obviously, we can assume that $h(p) > 0$). It follows that in condition VIII the constant $K$ is negative.

*Theorem* (on the evaluation of the residual term $R$).

If conditions I-VII are satisfied, there are constants $\varepsilon_0 > 0$ and $T_0 > 0, C \geq 0$, independent of $\varepsilon$, such that for all $0 < \varepsilon < \varepsilon_0$, $|x| < \infty$, $0 < t \leq T_0$ the solution of problem (4) exists, is unique, and the inequality holds:

$$\|R_N\| \leq C \varepsilon^{N+1}.$$

To prove the theorem, we construct the expansion $N+3$ term (which is guaranteed by sufficient smoothness of the initial condition), and represent the solution of problem (1) - (2) in the form:

$$U(x,t,p,\varepsilon) = S_{N+3}(\zeta,t,p) + P_{N+3}(\xi,\tau,p) + R =$$
$$= \sum_{i=0}^{N+3} \varepsilon^i (s_i(\zeta,t,p) + p_i(\xi,\tau,p)) + R = U_{N+3} + R_{N+3}, \quad (5)$$

where $U_{N+3}$ is the constructed asymptotics, $R_{N+3}$ is the residual term.

From the algorithm for constructing $U_{N+3}$ and evaluating functions $s_i, p_i$, it follows that the residual term $R_{N+3}$ satisfies the problem:

$$\varepsilon^2 (R_{N+3,t} + D(p)R_{N+3,x}) = L_p R_{N+3} + \varepsilon^2 RF(R_{N+3}, p) + \varepsilon^{N+3} f(x,t,p), \quad (6)$$

$$R_{N+3}(x,0,p) = 0; \quad (7)$$

$$f(x,t,p) = O(1). \quad (8)$$

The function $RF$ was defined as $RF(R,p) = F(S_{N+3} + P_{N+3} + R_{N+3}, p) - F(S_{N+3} + P_{N+3}, p)$.

Let us rewrite equation (6) as

$$\varepsilon^2 (R_t + DR_x) = r(x,t,p,R,\varepsilon). \quad (9)$$

Everywhere below, it is assumed that the variable $p$ changes on a certain segment $[p_1, p_2]$, and all functions are continuous with respect to this variable.

Let us take an arbitrary point $M_0(x_0, t_0, p_0)$ and omit from it the characteristics of equation (9) at different $p$ points on the axis $t = 0$. Let's denote $d_{max} = \max_p d(p), d_{min} = \min_p d(p)$. Obviously, the segment $[M_1 = x_0 - d_{min}t_0, M_2 = x_0 - d_{max}t_0]$ on the axis $t = 0$ will be the area of influence of the initial conditions, and the triangle with vertices at points $\{M_0, M_1, M_2\}$ will be the area of influence of the right side [6]. The resulting characteristic triangle is denoted by $\Delta_0$, and its boundary $t = 0$ is denoted by $\Gamma_0$.

The proof of the theorem is carried out by proving a number of lemmas.

*Lemma 1.*

Let $L$ be a linear invertible operator acting on a function $y = y(p)$, with the property: if the function $f(p)$ is positive $f_j(p) > 0, p_1 \le p \le p_2$, then the solution of the equation $Ly = f$ is positive: $y(p) > 0, p_1 \le p \le p_2$. Let $LY_1 = f_1, LY_2 = f_2$.

If $f_2(p) > |f_1(p)| \forall p,$ then $y_2(p) > |y_1(p)| \forall p.$

*Lemma 2.*

Let $U$ be the solution of the problem

$$u_t + Du_x = L_p u + f(x,t,p),$$
$$u(x,0,p) = u^0(x,p), \qquad (10)$$

where the operator $L$ satisfies condition VIII.

If for point $M_0(x_0, t_0, p_0)$ the inequalities for $f(x,t,p) > 0$ are satisfied in $\Delta_0$, $u_i^0(x,p) > 0$ on $\Gamma_0$ for all $p_1 \le p \le p_2$, then $u(x,t,p) > 0$ in $\Delta_0$ for all $p_1 \le p \le p_2$.

*Lemma 3* (on barriers)

Let $u_1, u_2$ be problem solving

$$(u_{it} + du_{ix}) = Lu_i + f_i(x,t,p), u_i(x,0,p) = u_i^0(x,p),$$

where $i = 1, 2$, the operator $L$ satisfies condition VIII.

If $u_1^0 > |u_2^0|, f_1 > |f_2|$, then $\forall (x,t,p): u_1 > |u_2|$.

*Lemma 4.*

Let $u$ be the solution of the problem (10), where the operator satisfies conditions VII, VIII. Then $\exists C > 0$, such that $\forall M_0(x_0, t_0, p_0)$ the inequality holds:

$$\|u\|_\Delta \le C(\|u^0\|_\Gamma + t_0 \|f\|_\Delta). \qquad (11)$$

*Lemma 5.*

Let $u$ be the solution of the problem

$$(u_t + du_x) = L_p u + f_1(x,t,p) + f_2(x,t,p,u), \qquad (12)$$
$$x_1, x_2, t \in \Omega = \{|x_1| < \infty, |x_2| < \infty, 0 < t < T\}, \qquad (13)$$
$$u(x,0,p) = u^0(x,p), \qquad (14)$$

where the operator $L_p$ satisfies conditions VII, VIII.

Let

1. The initial conditions and the function $f_1$ are continuous in their domains of definition.
2. Let the function $f_2(x,t,p,u)$ in the domain $G = \{|u| < K, K > 0\}$ be continuously differentiable with respect to the variable $u$ and $f_2(x,t,p,0) = 0$.
3. $|u^0| < C_1 K$, where $C_1 = 1/\min_p h(p))$.

Then there $\varepsilon_0 > 0$ also exist constants $0 < T_0 \le T, C > 0$ independent of $\varepsilon$, such that $\forall 0 < \varepsilon < \varepsilon_0, \forall M_0(x_{1_0}, x_{2_0}, t_0), 0 \le t_0 \le T_0$, the solution of problem (12) - (14) exists in $\Delta_0$, is unique and satisfies the bound

$$\|u\|_\Delta \le C\left(\|u^0\|_\Gamma + t_0 \|f_1\|_\Delta\right). \qquad (15)$$

*Remark.* The constant $C > 0$ in the evaluation (15) does not depend on the norm of the operator $L_p$.

Let $R_{N+3}$ be the solution of problem (6), (7). Let us divide equation (6) by $\varepsilon^2$

$$(R_{N+3,t} + D(p)R_{N+3,x}) = \varepsilon^{-2}L_p R_{N+3} + RF(R_{N+3}, p) + \varepsilon^{N+1}f(x,t,p),$$

Applying the estimate (15) of Lemma 5 to problem (6) - (7), taking into account the remark to Lemma 5, we obtain the estimate:

$$\|R_{N+3}\|_\Delta \leq C\left(\|R_{N+3}(x_1, x_2, 0)\| + t_0 \|\varepsilon^{N+1}f_1\|_\Delta\right);$$

Substituting $R_{N+3}(x_1, x_2, 0) = 0$ and $\|f_1\| = O(1)$, we get:

$$\|R_{N+3}\| = O(\varepsilon^{N+1}).$$

Because

$$R_N - R_{N+3} = \sum_{i=N+1}^{N+3} \varepsilon^i (s_i(\zeta, t, p) + p_i(\xi, \tau, p)) = O(\varepsilon^{N+1})$$

to i $\|R_{N+1}\| = \|R_{N+3} + O(\varepsilon^{N+1})\| = O(\varepsilon^{N+1})$

That's what we needed to prove.

## 3 Conclusion

The theorem on the uniform estimation of the residual term in AR (3) is proved. The next step is to extend the results to equation (1) with many spatial variables.